\documentclass[12pt]{article}
\usepackage{amssymb}
\usepackage{amsmath}
\usepackage{graphs}
\usepackage{amssymb}
\usepackage{amsthm}
\usepackage{url}
\usepackage[british]{babel}
\usepackage{bigstrut,array,float, fancyvrb,marvosym,verbatim,combgames}

\usepackage{psgo,array,url}

\theoremstyle{plain}
\newtheorem{definition}{Definition}
\newtheorem{lemma}[definition]{Lemma}
\newtheorem{theorem}[definition]{Theorem}

\newcommand{\lef}{\mathcal{L}}
\newcommand{\pre}{\mathcal{P}}
\newcommand{\n}{\mathcal{N}}
\newcommand{\ri}{\mathcal{R}}
\newcommand{\ti}{\mathcal{T}}

\newcommand{\gfr}{G_F^{SR}}
\newcommand{\gfl}{G_F^{SL}}
\newcommand{\plus}{+_{\ell}}

\newcommand{\tx}{\hbox{Tame}_x}
\newcommand{\tm}{\hbox{Tame}_{-x}}

\begin{document}
\title{Tame and Wild Scoring Play Games}
\author{Fraser Stewart\\
\small{Xi'An Jiaotong University, Department of Mathematics and Statistics, Xi'An, China}\\
\texttt{fraseridstewart@gmail.com}}
\date{}

\maketitle{}

\begin{abstract}
In this paper, we will be proving mathematically that scoring play combinatorial game theory covers all combinatorial games.  That is, there is a sub-set of scoring play games that are identical to the set of normal play games, and a different sub-set that is identical to the set of mis\`ere play games.

This proves conclusively, that scoring play combinatorial game theory is a complete theory for combinatorial games, and that every combinatorial game, regardless of the rule-set, can be analysed using scoring play combinatorial game theory.
\end{abstract}

\section{Introduction}

The theory of normal play games contains some of the most beautiful mathematics devised in the 20th century.  Out of it came the incredibly elegant ``surreal numbers'', as well as the foundation for the entire subject of combinatorial game theory.

Ever since Winning Ways \cite{ww} and On Numbers and Games \cite{onag} were published, mathematicians have believed that it is a complete theory for combinatorial games.  In this paper, we will show that the new theory of scoring play games as devised by Fraser Stewart \cite{fs}, can in fact be used to analyse all combinatorial games.

That is, there is a subset of scoring play games that is identical to the set of normal play games, and a further subset that is identical to the set of mis\`ere play games.  This means that scoring play theory is simply a more general version of the mathematics already described by Berlekamp, Conway and Guy.

\subsection{Scoring Play Combinatorial Game Theory}

One of the most recent major developments in the study of combinatorial games, was the paper ``Scoring Play Combinatorial Games'' \cite{fs}, which is currently due to appear in the volume Games of No Chance 5.  This paper offered a brand new mathematical theory for a class of games that had, until that point, received very little attention from combinatorial game theorists.

The first papers written on scoring play games were by Milnor and Hanner \cite{jm, oh}.  These were then followed up with papers with Ettinger \cite{je1, je2}, the first of which was published, while the second remains unpublished.  Lastly Will Johnson did some follow up work most recently in 2011 \cite{wj}.

However, the paper ``Scoring Play Combinatorial Games'', offers the most general theory for scoring play games.  All the classes of games studied by Milnor, Hanner, Ettinger and Johnson can be thought of as subclasses, of the class defined in this paper.  This work was done independently of the previous four authors, and is based on the theory of combinatorial games presented by Elwyn Berlekamp, John Conway and Richard Guy in Winning Ways and On Numbers and Games \cite{ww, onag}.

\begin{figure}[htb]
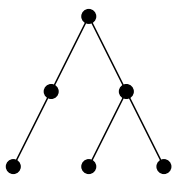


\begin{center}
\begin{graph}(2,2)

\roundnode{2}(0,0)\roundnode{3}(0.5,1)\roundnode{4}(1,2)
\roundnode{5}(1.5,1)\roundnode{6}(1,0)\roundnode{7}(2,0)

\edge{2}{3}\edge{3}{4}
\edge{4}{5}\edge{5}{6}\edge{5}{7}

\end{graph}
\end{center}

\caption{A typical game tree.}\label{gt}
\end{figure}

The idea behind this theory is very simple, consider the game tree given in figure \ref{gt}.

On a typical game tree, as shown in figure \ref{gt}, the nodes represent the positions of a game, and edges represent possible moves for both players from those positions.  Left sloping edges are Left's moves, and right sloping edges are Right's moves.

A scoring play game tree is exactly the same, but for one difference.  The nodes now have numbers on them which represent the score associated with that position.  The score is the difference between Left's total points, and Right's total points, at that point in the game.

\begin{figure}[htb]
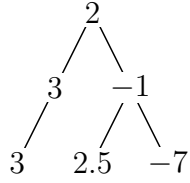


\begin{center}
\begin{graph}(2,2)

\roundnode{2}(0,0)\roundnode{3}(0.5,1)\roundnode{4}(1,2)
\roundnode{5}(1.5,1)\roundnode{6}(1,0)\roundnode{7}(2,0)

\edge{2}{3}\edge{3}{4}
\edge{4}{5}\edge{5}{6}\edge{5}{7}

\nodetext{2}{3}\nodetext{3}{3}\nodetext{4}{2}
\nodetext{5}{$-1$}\nodetext{6}{2.5}\nodetext{7}{$-7$}

\end{graph}
\end{center}

\caption{A scoring play game tree.}
\end{figure}

Formally scoring play games are defined as follows.

\begin{definition} A scoring play game $G=\{G^L|G^S|G^R\}$, where $G^L$ and $G^R$ are sets of games and $G^S\in\mathbb{R}$, the base case for the recursion is any game $G$ where $G^L=G^R=\emptyset$.\\

\noindent $G^L=\{\hbox{All games that Left can move to from } G\}$\\
\noindent $G^R=\{\hbox{All games that Right can move to from } G\}$,\\

and for all $G$ there is an $S=(P,Q)$ where $P$ and $Q$ are the number of points that Left and Right have on $G$ respectively.  Then $G^S=P-Q$, and for all $g^L\in G^L$, $g^R\in G^R$, there is a $p^L,p^R\in\mathbb{R}$ such that $g^{LS}=G^S+p^L$ and $g^{RS}=G^S+p^R$.

\end{definition}

By convention, we will take $G^S$ to be 0, unless stated otherwise.  This is simply to give games a ``default'' setting, i.e. if we don't know what $G^S$ is then it is natural to simply let it be 0.  We also write $\{.|G^S|.\}$ as $G^S$, e.g. $\{\{.|3|.\}|4|\{.|2|.\}\}$ would be written as $\{3|4|2\}$.  This simply for convenience and ease of reading.

For these games we also need the idea of a ``final score''.  That is the best possible score that both players can get when they move first.  Formally, this is defined as follows.

\begin{definition}  We define the following:

\begin{itemize}
\item{$\gfl$ is called the Left final score, and is the maximum score --when Left moves first on $G$-- at a terminal position on the game tree of $G$, if both Left and Right play perfectly.}
\item{$\gfr$ is called the Right final score, and is the minimum score --when Right moves first on $G$-- at a terminal position on the game tree of $G$, if both Left and Right play perfectly.}
\end{itemize}

\end{definition}

The paper ``Scoring Play Combinatorial Game Theory'' \cite{fs} discusses the structure of these games under the disjunctive sum, which is defined below.  In this paper it is shown that these games do not form a group, there is no non-trivial identity, and almost no games that can be compared in the usual sense.  However, these games are partially ordered under the disjunctive sum, and do form equivalence classes with a canonical form.  The games can also be reduced using the usual rules of domination and reversibility.

\begin{definition}  The disjunctive sum is defined as follows:
$$G\plus H=\{G^L\plus H,G\plus H^L|G^S+H^S|G^R\plus H,G\plus H^R\},$$
\noindent where $G^S+H^S$ is the normal addition of two real numbers.
\end{definition}

We abuse notation by letting $G^L$ and $G^R$ represent the set of options and the individual options themselves.  The reader will also notice that we have used $\plus$ and $+$.  This is to distinguish between the addition of games (the disjunctive sum), and the addition of scores.

\begin{definition} We define the following:
\begin{itemize}
\item{$-G=\{-G^R|-G^S|-G^L\}$.}
\item{For any two games $G$ and $H$, $G=H$ if $G\plus X$ has the same outcome as $H\plus X$ for all games $X$.}
\item{For any two games $G$ and $H$, $G\geq H$ if $H\plus X\in O$ implies $G\plus X\in O$, where $O=L_\geq$, $R_\geq$, $L_>$ or $R_>$, for all games $X$.}
\item{For any two games $G$ and $H$, $G\leq H$ if $H\plus X\in O$ implies $G\plus X\in O$, where $O=L_\leq$, $R_\leq$, $L_<$ or $R_<$, for all games $X$.}
\item{$G\cong H$ means $G$ and $H$ have identical game trees.}
\item{$G\approx H$ means $G$ and $H$ have the same outcome.}
\end{itemize}
\end{definition}  

We also define the outcome classes a little differently, we first need the following definition.

\begin{definition}
\item{$L_>=\{G|\gfl>0\}$, $L_<=\{G|\gfl<0\}$, $L_= =\{G|\gfl=0\}$.}
\item{$R_>=\{G|\gfr>0\}$, $R_<=\{G|\gfr<0\}$, $R_= =\{G|\gfr=0\}$.}
\item{$L_\geq=L_>\cup L_=$, $L_\leq = L_<\cup L_=$.}
\item{$R_\geq=R_>\cup R_=$, $L_\leq = R_<\cup R_=$.}
\end{definition}

This allows us to define the outcome classes as follows. 

\begin{definition}
The outcome classes of scoring games are defined as follows:

\begin{itemize}
\item{$\lef=(L_>\cap R_>)\cup(L_>\cap R_=)\cup(L_=\cap R_>)$}
\item{$\ri=(L_<\cap R_<)\cup(L_<\cap R_=)\cup(L_=\cap R_<)$}
\item{$\n=L_>\cap R_<$}
\item{$\pre=L_<\cap R_>$}
\item{$\ti =L_=\cap R_=$}
\end{itemize}
\end{definition}

Note that under these definitions, every game belongs to exactly one outcome class.

\section{Tame and Wild Games}

The terms ``tame'' and ``wild'' are taken directly from the Winning Ways \cite{ww}.  The authors in that book used it in the context of mis\`ere impartial games, in the sense that ``tame'' meant ``like nim'', and wild meant ``unlike nim''.

Here we are using the terms in the context of scoring play games to mean ``like normal play'', and ``unlike normal play''.  Those are very vague terms, but we will be giving a formal mathematical definition of tame and wild shortly.

To motivate this, first consider the game Pirates and Treasure, as first introduced in \cite{fs3}.  The rules are as follows:

\begin{enumerate}
\item{The game is played on a finite simple graph, defined arbitrarily before the game begins.  Left has $n$ ships, and Right has $m$ ship.}
\item{Each ship has a pre-defined starting vertex.}
\item{Every node is numbered to indicate how much treasure there is at that node, the players starting vertices are not numbered.}
\item{On a player's turn he moves to an adjacent, unvisited vertex.  The number of points he gets, corresponds to the number on that vertex. A player may not move a previously visited vertex, including the starting vertices.}
\item{The game ends when it is a player's turn and he is not adjacent to an unvisited vertex.}
\item{The player who gathers the most treasure wins.}
\end{enumerate}

In his paper, the author first showed that this game, in general, has no nice mathematical structure.  He then showed, that by simply fixing the values of the vertices to be some number $x>0$, this game starts to behave very similarly to a normal play game.

What we will do in this paper, is consider this class of games in general.  That is scoring games, which behave like a normal play, i.e. moving last is the player's best strategy.  First we give a formal definition.

\begin{definition}
A scoring play combinatorial game $G$ is tame if the score is greater than or equal to 0, whenever Left moves last, and less than or equal to 0, whenever Right moves last.
\end{definition}

This definition is a little stronger than simply saying $G^{SL}_F\geq 0$, if Left moves last, and likewise for the remaining cases.  In our definition we are stating that this is true for the entire game tree, i.e. no matter which options Left chooses, if he moves last then he will at least tie the game.  So if $G$ is tame then this implies that $G^{SL}_F\geq 0$ if Left moves last, but the converse is not true, and likewise for the remaining cases.

First, we will look at the general mathematical structure of this set, and ask if it as nice as the set of normal play games.  In other words, does making it advantageous to move last really affect the mathematical structure of scoring games.

\begin{theorem}\label{t1}  The outcome class table for tame games $G$, under the disjunctive sum, is given in Table \ref{table1}.

\begin{table}[htb]
\begin{center}
\begin{tabular}{c|ccccc}
$G\plus H$&$G\in\pre$&$G\in\ti$&$G\in\lef$&$G\in\ri$&$G\in\n$\\\hline
$H\in\pre$&$\pre$&$\pre,\lef,\ri,\n$&$\lef$&$\ri$&$\n$\\
$H\in\ti$&$\pre,\lef, \ri, \n$&$\pre,\ti,\lef,\ri,\n$&$\lef,\ti$&$\ri,\ti$&$\pre,\ti, \lef, \ri, \n$\\
$H\in\lef$&$\lef$&$\lef,\ti$&$\lef$&$\pre,\ti,\lef,\ri,\n$&$\lef,\n$\\
$H\in\ri$&$\ri$&$\ri,\ti$&$\pre,\ti,\lef,\ri,\n$&$\ri$&$\ri,\n$\\
$H\in\n$&$\n$&$\pre,\ti, \lef, \ri, \n$&$\lef, \n$&$\ri,\n$&$\pre,\ti,\lef,\ri,\n$\\
\end{tabular}
\end{center}
\caption{The outcome of $G\plus H$, where $G$ and $H$ are tame games.}\label{table1}
\end{table}

\end{theorem}

\begin{proof}
Case 1: $G,H\in\lef$ implies that $G\plus H\in\lef$.

From the definition of the outcome classes, we know that if $G$ and $H$ are both in $\lef$, then this implies that Left has a winning strategy either moving first or second, on both $G$ and $H$.

If Left has a winning strategy moving first on $G$, say, then Left can win moving first on $G\plus H$ by simply playing his winning move on $G$, i.e. moving to $G^L\plus H$.  He then always moves in the same component as Right, which will guarantee that he moves last on both $G$ and $H$.  Since the score of $G$ will be greater than zero, and the score of $H$ greater than, or equal, to zero, then the score of $G\plus H$ will be greater than zero.  Therefore, if Left can win moving first on $G$ or $H$, Left can win moving first on $G\plus H$.

If Left has a winning strategy moving second on $G$, but not $H$ (i.e. $H$ ends in a tie), then Left simply moves on the same component that Right moves on.  Since Left moves last on $H$, this will force Right to move first on $G$, which Left wins.  Therefore, if Left has a winning strategy moving second on either $G$ or $H$, then Left has a winning strategy moving second on $G\plus H$.

Since Left can guarantee a tie moving first or second, by an identical argument as above.  Then this implies that $G\plus H\in\lef$.

Case 2: $G\in\pre$, $H\in\mathcal{X}$ implies that $G\plus H\in\mathcal{X}$, where $\mathcal{X}=\lef, \ri, \pre$ and $\n$.

Consider Left moving first, since the case Right moving first will follow by symmetry.  If Left moves last on $H$ moving first, then Left will play to $G\plus H^L$.  Left then simply plays on the same component as Right, and he is guaranteed to move last on both $G$ and $H$, and since $G\in\pre$, the final score of $G\plus H$ will be greater than the final score of $H$.

If Right moves last on $H$ with Left moving first, then again Right can also move last on $G\plus H$, by simply moving on the same component as Left.  As before, since $G\in\pre$, the final score of $G\plus H$ will be lower than the final score of $H$.

Therefore, the outcome of $G\plus H$ will be the same as the outcome of $H$ if $H\in\lef,\ri,\pre$ or $\n$.

Case 3: $G\in\pre$, $H\in\ti$ implies that $G\plus H\in\mathcal{X}$, where $\mathcal{X}=\lef, \ri, \pre$ and $\n$.

The first thing is to show that if $G\in\pre$ and $H\in\ti$, then $G\plus H\not\in\ti$.  To prove this, consider Left moving first, since Right moving first follows by symmetry.  If Left moving first moves last  then Left moving first on $G\plus H$ will move to $G\plus H^L$, since he can simply move on the same component as Right and win outright.

If Left does not move last and moves to $G^L\plus H$, then Right will again simply respond on whichever component Left played on and wins.  If Left moves to $G\plus H^L$, then Right does the same thing.  In other words, when Left moves last on $H$, Right has no strategy to tie $G\plus H$, and likewise when Right moves last on $H$, Left has no strategy to tie $G\plus H$.  Therefore $G\plus H\not\in \ti$.  To complete this case, we give examples to show that there is a $G\in \pre, H\in\ti$ such that $G\plus H\in\mathcal{X}$, where $\mathcal{X}=\lef, \ri, \pre$ and $\n$.

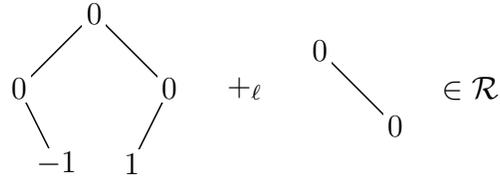
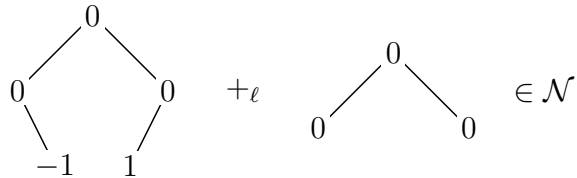
\begin{center}
\begin{graph}(6,2)

\roundnode{1}(4,1.5)\roundnode{2}(5,0.5)\roundnode{3}(1,2)\roundnode{4}(0,1)
\roundnode{5}(2,1)\roundnode{6}(0.5,0)\roundnode{7}(1.5,0)

\edge{1}{2}\edge{3}{4}\edge{3}{5}\edge{4}{6}\edge{5}{7}

\nodetext{1}{$0$}\nodetext{2}{0}\nodetext{3}{0}\nodetext{4}{0}
\nodetext{5}{0}\nodetext{6}{$-1$}\nodetext{7}{$1$}

\freetext(3,1){$\plus$}\freetext(6,1){$\in\ri$}

\end{graph}
\end{center}

\begin{center}
\begin{graph}(7,3)

\roundnode{1}(5,1.5)\roundnode{2}(6,0.5)\roundnode{3}(1,2)\roundnode{4}(0,1)
\roundnode{5}(2,1)\roundnode{6}(0.5,0)\roundnode{7}(1.5,0)\roundnode{8}(4,0.5)

\edge{1}{2}\edge{3}{4}\edge{3}{5}\edge{4}{6}\edge{5}{7}\edge{1}{8}

\nodetext{1}{$0$}\nodetext{2}{0}\nodetext{3}{0}\nodetext{4}{0}
\nodetext{5}{0}\nodetext{6}{$-1$}\nodetext{7}{$1$}\nodetext{8}{0}

\freetext(3,1){$\plus$}\freetext(7,1){$\in\n$}

\end{graph}
\end{center}

The case $G\in\pre$, $H\in\ti$ such that $G\plus H\in\pre$ is trivial, since we simply let $H=\{.|0|.\}$.  The remaining case follows by symmetry.

\noindent Case 4: $G\in\ti$, $H\in\lef$ implies that $G\plus H\in\pre, \ti, \lef,\ri$ or $\n$.

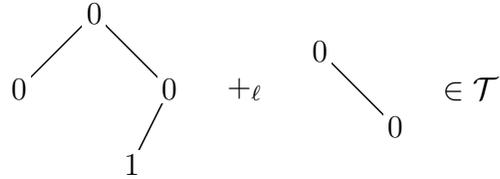
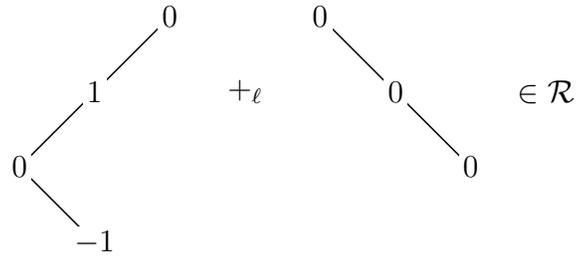
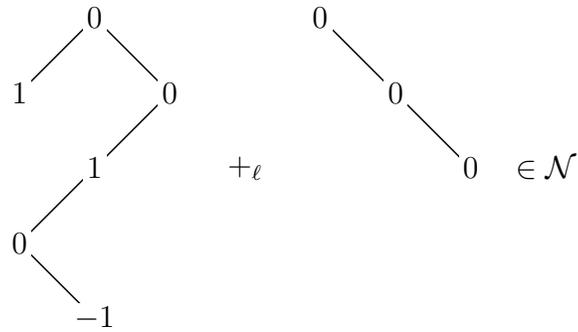
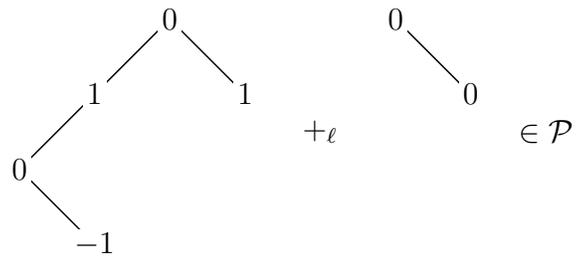
\begin{center}
\begin{graph}(6,2)

\roundnode{1}(1,2)\roundnode{2}(2,1)\roundnode{3}(4,1.5)\roundnode{4}(0,1)
\roundnode{5}(5,0.5)\roundnode{7}(1.5,0)

\edge{1}{2}\edge{1}{4}\edge{3}{5}\edge{2}{7}

\nodetext{1}{$0$}\nodetext{2}{0}\nodetext{3}{0}\nodetext{4}{0}
\nodetext{5}{0}\nodetext{7}{$1$}

\freetext(3,1){$\plus$}\freetext(6,1){$\in\ti$}

\end{graph}
\end{center}

\begin{center}
\begin{graph}(7,4)

\roundnode{1}(0,1)\roundnode{2}(1,0)\roundnode{3}(1,2)\roundnode{4}(2,3)
\roundnode{5}(4,3)\roundnode{6}(5,2)\roundnode{7}(6,1)

\edge{1}{2}\edge{1}{3}\edge{3}{4}\edge{5}{6}\edge{6}{7}

\nodetext{1}{0}\nodetext{2}{$-1$}\nodetext{3}{1}\nodetext{4}{0}
\nodetext{5}{0}\nodetext{6}{0}\nodetext{7}{0}

\freetext(3,2){$\plus$}\freetext(7,2){$\in\ri$}

\end{graph}
\end{center}

\begin{center}
\begin{graph}(7,5)

\roundnode{1}(0,3)\roundnode{2}(1,4)\roundnode{3}(2,3)\roundnode{4}(1,2)
\roundnode{5}(0,1)\roundnode{6}(1,0)\roundnode{7}(4,4)\roundnode{8}(5,3)
\roundnode{9}(6,2)

\edge{1}{2}\edge{2}{3}\edge{3}{4}\edge{4}{5}
\edge{5}{6}\edge{7}{8}\edge{8}{9}

\nodetext{1}{1}\nodetext{2}{0}\nodetext{3}{0}\nodetext{4}{1}
\nodetext{5}{0}\nodetext{6}{$-1$}\nodetext{7}{0}\nodetext{8}{0}
\nodetext{9}{0}

\freetext(3,2){$\plus$}\freetext(7,2){$\in\n$}

\end{graph}
\end{center}

\begin{center}
\begin{graph}(7,4)

\roundnode{1}(1,0)\roundnode{2}(0,1)\roundnode{3}(1,2)\roundnode{4}(2,3)
\roundnode{5}(3,2)\roundnode{6}(5,3)\roundnode{7}(6,2)

\edge{1}{2}\edge{2}{3}\edge{3}{4}\edge{4}{5}
\edge{6}{7}

\nodetext{1}{$-1$}\nodetext{2}{0}\nodetext{3}{1}\nodetext{4}{0}
\nodetext{5}{1}\nodetext{6}{0}\nodetext{7}{0}

\freetext(4,1.5){$\plus$}\freetext(7,1.5){$\in\pre$}

\end{graph}
\end{center}

\noindent Case 5: $G\in\ti$, $H\in\n$ implies that $G\plus H\in\mathcal{X}$, where $\mathcal{X}=\pre,\ti,\lef,\ri$ or $\n$.

\begin{center}
\begin{graph}(8,3)

\roundnode{1}(0,1.5)\roundnode{2}(1,2.5)\roundnode{3}(2,1.5)\roundnode{4}(4,1)
\roundnode{5}(5,0)\roundnode{6}(5,2)\roundnode{7}(6,3)\roundnode{8}(7,2)

\edge{1}{2}\edge{2}{3}\edge{4}{5}
\edge{4}{6}\edge{6}{7}\edge{7}{8}

\nodetext{1}{0}\nodetext{2}{0}\nodetext{3}{0}\nodetext{4}{0}
\nodetext{5}{$-1$}\nodetext{6}{1}\nodetext{7}{0}\nodetext{8}{$-1$}

\freetext(3,1.5){$\plus$}\freetext(8,1.5){$\in\ri$}

\end{graph}
\end{center}

\begin{center}
\begin{graph}(9,4)

\roundnode{1}(0,1.5)\roundnode{2}(1,2.5)\roundnode{3}(2,1.5)\roundnode{4}(4,1)
\roundnode{5}(5,0)\roundnode{6}(5,2)\roundnode{7}(6,3)\roundnode{8}(7,2)\roundnode{9}(8,1)
\roundnode{10}(7,0)

\edge{1}{2}\edge{2}{3}\edge{4}{5}
\edge{4}{6}\edge{6}{7}\edge{7}{8}
\edge{8}{9}\edge{9}{10}

\nodetext{1}{0}\nodetext{2}{0}\nodetext{3}{0}\nodetext{4}{0}
\nodetext{5}{$-1$}\nodetext{6}{1}\nodetext{7}{0}\nodetext{8}{$-1$}
\nodetext{9}{0}\nodetext{10}{1}

\freetext(3,1.5){$\plus$}\freetext(9,1.5){$\in\pre$}

\end{graph}
\end{center}

\begin{center}
\begin{graph}(9,3)

\roundnode{1}(0,0.5)\roundnode{2}(1,1.5)\roundnode{3}(2,0.5)\roundnode{4}(4,0)
\roundnode{6}(5,1)\roundnode{7}(6,2)\roundnode{8}(7,1)\roundnode{9}(8,0)

\edge{1}{2}\edge{2}{3}
\edge{4}{6}\edge{6}{7}\edge{7}{8}
\edge{8}{9}

\nodetext{1}{0}\nodetext{2}{0}\nodetext{3}{0}\nodetext{4}{0}
\nodetext{6}{1}\nodetext{7}{0}\nodetext{8}{$-1$}
\nodetext{9}{0}

\freetext(3,0.5){$\plus$}\freetext(9,0.5){$\in\ti$}

\end{graph}
\end{center}

\noindent Case 6: $G\in\lef$, $H\in\ri$ implies that $G\plus H\in \mathcal{X}$, where $\mathcal{X}=\pre, \ti, \lef, \ri$ or $\n$.

\begin{center}
\begin{graph}(6,2)

\roundnode{1}(1,2)\roundnode{2}(2,1)\roundnode{3}(4,1.5)\roundnode{4}(0,1)
\roundnode{5}(5,0.5)\roundnode{7}(1.5,0)

\edge{1}{2}\edge{1}{4}\edge{3}{5}\edge{2}{7}

\nodetext{1}{$0$}\nodetext{2}{0}\nodetext{3}{0}\nodetext{4}{1}
\nodetext{5}{$-1$}\nodetext{7}{$1$}

\freetext(3,1){$\plus$}\freetext(6,1){$\in\ti$}

\end{graph}
\end{center}

\begin{center}
\begin{graph}(7,4)

\roundnode{1}(0,1)\roundnode{2}(1,0)\roundnode{3}(1,2)\roundnode{4}(2,3)
\roundnode{5}(4,3)\roundnode{6}(5,2)\roundnode{7}(6,1)

\edge{1}{2}\edge{1}{3}\edge{3}{4}\edge{5}{6}\edge{6}{7}

\nodetext{1}{0}\nodetext{2}{$-1$}\nodetext{3}{1}\nodetext{4}{0}
\nodetext{5}{0}\nodetext{6}{$-1$}\nodetext{7}{$-1$}

\freetext(3,2){$\plus$}\freetext(7,2){$\in\ri$}

\end{graph}
\end{center}

\begin{center}
\begin{graph}(7,5)

\roundnode{1}(0,3)\roundnode{2}(1,4)\roundnode{3}(2,3)\roundnode{4}(1,2)
\roundnode{5}(0,1)\roundnode{6}(1,0)\roundnode{7}(4,4)\roundnode{8}(5,3)
\roundnode{9}(6,2)

\edge{1}{2}\edge{2}{3}\edge{3}{4}\edge{4}{5}
\edge{5}{6}\edge{7}{8}\edge{8}{9}

\nodetext{1}{1}\nodetext{2}{0}\nodetext{3}{0}\nodetext{4}{1}
\nodetext{5}{0}\nodetext{6}{$-1$}\nodetext{7}{$-1$}\nodetext{8}{0}
\nodetext{9}{0}

\freetext(3,2){$\plus$}\freetext(7,2){$\in\n$}

\end{graph}
\end{center}

\begin{center}
\begin{graph}(7,4)

\roundnode{1}(1,0)\roundnode{2}(0,1)\roundnode{3}(1,2)\roundnode{4}(2,3)
\roundnode{5}(3,2)\roundnode{6}(5,3)\roundnode{7}(6,2)

\edge{1}{2}\edge{2}{3}\edge{3}{4}\edge{4}{5}
\edge{6}{7}

\nodetext{1}{$-1$}\nodetext{2}{0}\nodetext{3}{1}\nodetext{4}{0}
\nodetext{5}{1}\nodetext{6}{$-1$}\nodetext{7}{0}

\freetext(4,1.5){$\plus$}\freetext(7,1.5){$\in\pre$}

\end{graph}
\end{center}

\noindent Case 7: $G\in\lef$, $H\in\n$ implies that $G\plus H\in \lef$ or $\n$.

Left can win moving first on $G\plus H$ by moving to $G\plus H^L$, since he simply moves on the same component as Right and can move last on $G$ and $H$.  Since $H\in\n$ the final score of $H$ will be greater than 0, and the final score of $G$ is greater than or equal to zero, then the final score of $G\plus H$ will be greater than 0.

Since Left always wins moving first on $G\plus H$, we can conclude that $G\plus H\not\in\pre,\ti$ or $\ri$.  To complete the proof we give an example of $G\plus H\in\lef$ and $G\plus H\in\n$.

\begin{center}
\begin{graph}(7,3)

\roundnode{1}(0,1)\roundnode{2}(1,2)\roundnode{3}(2,1)\roundnode{4}(1.5,0)
\roundnode{5}(4,1)\roundnode{6}(5,2)\roundnode{7}(6,1)

\edge{1}{2}\edge{2}{3}\edge{3}{4}\edge{5}{6}\edge{6}{7}

\nodetext{1}{1}\nodetext{2}{0}\nodetext{3}{0}\nodetext{4}{1}
\nodetext{5}{1}\nodetext{6}{0}\nodetext{7}{$-1$}

\freetext(3,1){$\plus$}\freetext(7,1){$\in\lef$}

\end{graph}
\end{center}

\begin{center}
\begin{graph}(7,3)

\roundnode{1}(0,1)\roundnode{2}(1,2)\roundnode{3}(2,1)\roundnode{4}(1.5,0)
\roundnode{5}(4,1)\roundnode{6}(5,2)\roundnode{7}(6,1)

\edge{1}{2}\edge{2}{3}\edge{3}{4}\edge{5}{6}\edge{6}{7}

\nodetext{1}{1}\nodetext{2}{0}\nodetext{3}{0}\nodetext{4}{1}
\nodetext{5}{1}\nodetext{6}{0}\nodetext{7}{$-2$}

\freetext(3,1){$\plus$}\freetext(7,1){$\in\n$}

\end{graph}
\end{center}

\noindent Case 8: $G, H\in\n$ implies that $G\plus H\in \mathcal{X}$, where $\mathcal{X}=\pre, \ti, \lef, \ri$ or $\n$.

\begin{center}
\begin{graph}(7,2)

\roundnode{1}(0,0)\roundnode{2}(1,1)\roundnode{3}(2,0)
\roundnode{4}(4,0)\roundnode{5}(5,1)\roundnode{6}(6,0)

\edge{1}{2}\edge{2}{3}
\edge{4}{5}\edge{5}{6}

\nodetext{1}{2}\nodetext{2}{0}\nodetext{3}{$-1$}
\nodetext{4}{2}\nodetext{5}{0}\nodetext{6}{$-1$}

\freetext(3,0.5){$\plus$}\freetext(7,0.5){$\in\lef$}

\end{graph}
\end{center}

\begin{center}
\begin{graph}(7,2)

\roundnode{1}(0,0)\roundnode{2}(1,1)\roundnode{3}(2,0)
\roundnode{4}(4,0)\roundnode{5}(5,1)\roundnode{6}(6,0)

\edge{1}{2}\edge{2}{3}
\edge{4}{5}\edge{5}{6}

\nodetext{1}{1}\nodetext{2}{0}\nodetext{3}{$-1$}
\nodetext{4}{1}\nodetext{5}{0}\nodetext{6}{$-1$}

\freetext(3,0.5){$\plus$}\freetext(7,0.5){$\in\ti$}

\end{graph}
\end{center}

\begin{center}
\begin{graph}(7,2)

\roundnode{1}(0,0)\roundnode{2}(1,1)\roundnode{3}(2,0)
\roundnode{4}(4,0)\roundnode{5}(5,1)\roundnode{6}(6,0)

\edge{1}{2}\edge{2}{3}
\edge{4}{5}\edge{5}{6}

\nodetext{1}{2}\nodetext{2}{0}\nodetext{3}{$-2$}
\nodetext{4}{1}\nodetext{5}{0}\nodetext{6}{$-1$}

\freetext(3,0.5){$\plus$}\freetext(7,0.5){$\in\n$}

\end{graph}
\end{center}

\begin{center}
\begin{graph}(9,4)

\roundnode{1}(0,1.5)\roundnode{2}(1,2.5)\roundnode{3}(2,1.5)\roundnode{4}(4,1)
\roundnode{5}(5,0)\roundnode{6}(5,2)\roundnode{7}(6,3)\roundnode{8}(7,2)\roundnode{9}(8,1)
\roundnode{10}(7,0)

\edge{1}{2}\edge{2}{3}\edge{4}{5}
\edge{4}{6}\edge{6}{7}\edge{7}{8}
\edge{8}{9}\edge{9}{10}

\nodetext{1}{$\frac{1}{2}$}\nodetext{2}{0}\nodetext{3}{$-\frac{1}{2}$}\nodetext{4}{0}
\nodetext{5}{$-1$}\nodetext{6}{1}\nodetext{7}{0}\nodetext{8}{$-1$}
\nodetext{9}{0}\nodetext{10}{1}

\freetext(3,1.5){$\plus$}\freetext(9,1.5){$\in\pre$}

\end{graph}
\end{center}

All remaining cases follow by symmetry, and so this completes the proof.

\end{proof}

\begin{theorem}
Tame scoring play games form a non-trivial monoid under the disjunctive sum.
\end{theorem}

\begin{proof}  

\noindent Part 1: Identity

If $X\in \ti$ and the second player moves last on $X$, then $G\plus X\approx G$ for all tame games $G$.

There are four cases to consider, since the other four follow by symmetry, Left wins moving first on $G$, ties moving first on $G$, ties moving second on $G$, or loses moving first on $G$.  When we say a player ties moving first or second, this means that the player is the last player to move, and the score was a tie.

If Left wins moving first on $G$, then when Left plays $G\plus X$, Left can win by simply moving to $G^L\plus X$.  Since Left moves last on both $G$ and $X$, then whichever component Right chooses, Left plays the same component and is guaranteed to win.

If Left ties moving first on $G$, then again Left can tie playing exactly the same strategy as above.  Left cannot win either $G$ or $X$, therefore playing this strategy guarantees a tie, and so $G\plus X$ will also finish in a tie for Left.

If Left ties moving second on $G$, then again by the same argument, Left will play to move last on $G\plus X$.  Left can do so simply by moving on the same component as Right.  By the same argument as above, Right will first move to $G^R\plus X$, as this is his best move.

If Left loses moving first, then by the same argument, Right can guarantee a win by simply playing the same component as Left at each turn.

\noindent Part 2: Inverse

We take $-G$, to be the inverse of $G$.  To show that this is the inverse, we will show that $G\plus (-G)\plus X\approx X$ for all tame games $X$.  The argument for this is almost identical to Part 1.  The difference being that if, Left say, chooses to move to $G^L\plus (-G)\plus X$, Right will move to $G^L\plus (-G^R)\plus X$.  Apart from that, the argument is completely identical.

\noindent Part 3: Not Closed.

To that this set is a monoid, and not a group, we need to show that it is not closed.  To do this consider the game $G=\{a|b|.\}$ and $H=\{.|c|d\}$.  By definition $a,b\geq 0$ and $c,d\leq 0$, so simply choose $a>|b|$, then Left moving first on $G\plus H$ moves to $a\plus H$, Right must move to $a\plus d$.  But, $a+d>0$, and Right moved last on $G\plus H$.  Therefore, $G\plus H$ is not a tame game, and the set is not closed under the disjunctive sum.

So the theorem is proven.

\end{proof}

\subsection{Mapping from Normal Play}\label{normal}

The real question is ``Can we use normal play theory on these games?''.  What this question is really asking is ``Is there a mapping from the set of normal play games, to the set of scoring play games?''.  In other words, winning under normal play, implies winning under scoring play.  We begin by defining the following set.

\begin{definition}\label{tame}
$\hbox{Tame}_x=\{G|G=\{G^L\plus x|0|G^R\plus -x\}, G^L,G^R\in\hbox{Tame}_x\}$
\end{definition}

The aim here is to show that under the definitions of scoring play theory, the set $\hbox{Tame}_x$ is in fact, completely equivalent, to the set of normal play combinatorial games.  That is, the two sets are identical.  Further, we can use definition \ref{tame} to derive surreal numbers.

\begin{theorem}\label{t2} The outcome class table for games in $\hbox{Tame}_x$ is given in Table \ref{table2}.

\begin{table}[htb]
\begin{center}
\begin{tabular}{c|cccc}
$G\plus H$&$G\in\ti$&$G\in\lef$&$G\in\ri$&$G\in\n$\\\hline
$H\in\ti$&$\ti$&$\lef$&$\ri$&$\n$\\
$H\in\lef$&$\lef$&$\lef$&$\lef$, $\ri$, $\n$, $\ti$&$\lef$, $\n$\\
$H\in\ri$&$\ri$&$\lef$, $\ri$, $\n$, $\ti$&$\ri$&$\ri$, $\n$\\
$H\in\n$&$\n$&$\lef$, $\n$&$\ri$, $\n$&$\lef$, $\ri$,$\n$, $\ti$\\
\end{tabular}
\end{center}
\caption{Outcome Class Table for $\hbox{Tame}_x$}\label{table2}
\end{table}

\end{theorem}

\begin{proof} 

Case 1:  $G\in\ti$, $H\in\mathcal{X}$ implies that $G\plus H\in\mathcal{X}$, where $\mathcal{X}=\ti,\lef,\ri$ or $\n$.

Consider the game $G\plus H$, we only discuss the cases where Left moves last moving first and second on $H$, since the other cases follow by symmetry.  Left moving first and moving last on $H$ implies that the score of $H$ must be $x$ when it ends.

So when Left plays $G\plus H$, then he can win by moving to $G\plus H^L$.  Since he can move last on both $G$ and $H$, he simply plays on whichever component Right plays on.  Therefore the final score of $G\plus H$ will also be $x$, when Left plays first.

If Left moves last moving second, then Right's best move, is to move to $G\plus H^R$.  Since final score of $G\plus H$ will either be $-x$ or $0$, when Right moves first, then Right cannot win by moving to $G^R\plus H$.  If Right does this, then Left can move to $G^R\plus H^L$ and the final score of this game will be greater than or equal to zero.

If Right moves to $G\plus H^R$, then he can guarantee a tie by simply moving on the same component as Left.  Therefore $G\plus X\approx G$ for all games $G\in\hbox{Tame}_x$.

\noindent Case 2: $G,H\in\lef$ implies that $G\plus H\in\lef$.

This follows from theorem \ref{t1}.

\noindent Case 3: $G\in\lef$, $H\in\n$ implies that $G\plus H\in\lef$ or $\n$.

This also follows from theorem \ref{t1}.

\noindent Case 4: $G,H\in\n$ implies that $G\plus H\in\mathcal{X}$, where $\mathcal{X}=\ti,\lef,\ri$ or $\n$.

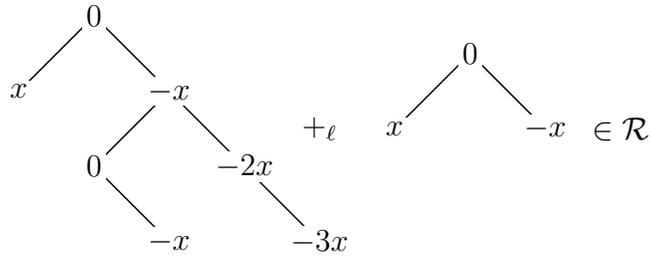
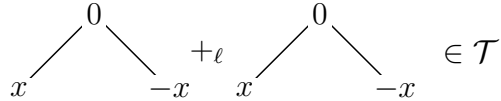
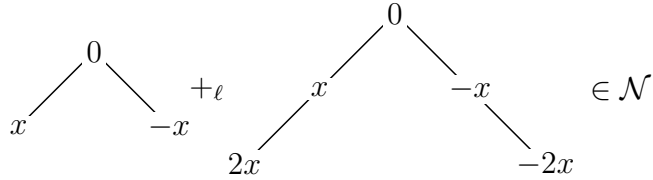
\begin{center}
\begin{graph}(8,4)

\roundnode{1}(0,2)\roundnode{2}(1,3)\roundnode{3}(2,2)\roundnode{4}(1,1)\roundnode{5}(2,0)
\roundnode{6}(3,1)\roundnode{7}(4,0)\roundnode{8}(5,1.5)\roundnode{9}(6,2.5)\roundnode{10}(7,1.5)

\edge{1}{2}\edge{2}{3}\edge{3}{4}\edge{4}{5}
\edge{3}{6}\edge{6}{7}\edge{8}{9}\edge{9}{10}

\nodetext{1}{$x$}\nodetext{2}{0}\nodetext{3}{$-x$}\nodetext{4}{0}\nodetext{5}{$-x$}
\nodetext{6}{$-2x$}\nodetext{7}{$-3x$}\nodetext{8}{$x$}\nodetext{9}{0}\nodetext{10}{$-x$}

\freetext(4,1.5){$\plus$}\freetext(8,1.5){$\in\ri$}

\end{graph}
\end{center}

\begin{center}
\begin{graph}(6,2)

\roundnode{1}(0,0)\roundnode{2}(1,1)\roundnode{3}(2,0)
\roundnode{4}(3,0)\roundnode{5}(4,1)\roundnode{6}(5,0)

\edge{1}{2}\edge{2}{3}\edge{4}{5}\edge{5}{6}

\nodetext{1}{$x$}\nodetext{2}{0}\nodetext{3}{$-x$}
\nodetext{4}{$x$}\nodetext{5}{0}\nodetext{6}{$-x$}

\freetext(2.5,0.5){$\plus$}\freetext(6,0.5){$\in\ti$}

\end{graph}
\end{center}

\begin{center}
\begin{graph}(8,3)

\roundnode{1}(0,0.5)\roundnode{2}(1,1.5)\roundnode{3}(2,0.5)\roundnode{4}(3,0)
\roundnode{5}(4,1)\roundnode{6}(5,2)\roundnode{7}(6,1)\roundnode{8}(7,0)

\edge{1}{2}\edge{2}{3}\edge{4}{5}
\edge{5}{6}\edge{6}{7}\edge{7}{8}

\nodetext{1}{$x$}\nodetext{2}{0}\nodetext{3}{$-x$}\nodetext{4}{$2x$}
\nodetext{5}{$x$}\nodetext{6}{0}\nodetext{7}{$-x$}\nodetext{8}{$-2x$}

\freetext(2.5,1){$\plus$}\freetext(8,1){$\in\n$}

\end{graph}
\end{center}

\noindent Case 5: $G\in\lef$, $H\in\ri$ implies that $G\plus H\in\mathcal{X}$, where $\mathcal{X}=\ti,\lef,\ri$ or $\n$.

\begin{center}
\begin{graph}(7,3)

\roundnode{1}(0,1)\roundnode{2}(1,2)\roundnode{3}(2,1)\roundnode{4}(1.5,0)
\roundnode{5}(4.5,0)\roundnode{6}(4,1)\roundnode{7}(5,2)\roundnode{8}(6,1)

\edge{1}{2}\edge{2}{3}\edge{3}{4}
\edge{5}{6}\edge{6}{7}\edge{7}{8}

\nodetext{1}{$x$}\nodetext{2}{0}\nodetext{3}{$-x$}\nodetext{4}{0}
\nodetext{5}{0}\nodetext{6}{$x$}\nodetext{7}{0}\nodetext{8}{$-x$}

\freetext(3,1){$\plus$}\freetext(7,1){$\in\ti$}

\end{graph}
\end{center}

\begin{center}
\begin{graph}(6,3)

\roundnode{1}(0,1)\roundnode{2}(1,2)\roundnode{3}(3,2)\roundnode{4}(4,1)\roundnode{5}(5,0)

\edge{1}{2}\edge{3}{4}\edge{4}{5}

\nodetext{1}{0}\nodetext{2}{$x$}\nodetext{3}{0}\nodetext{4}{$-x$}\nodetext{5}{$-2x$}

\freetext(2,1){$\plus$}\freetext(6,1){$\in\ri$}

\end{graph}
\end{center}

\begin{center}
\begin{graph}(6,3)

\roundnode{1}(0,0)\roundnode{2}(1,1)\roundnode{3}(2,2)\roundnode{4}(4,2)\roundnode{5}(5,1)

\edge{1}{2}\edge{2}{3}\edge{4}{5}

\nodetext{1}{$2x$}\nodetext{2}{$x$}\nodetext{3}{0}\nodetext{4}{0}\nodetext{5}{$-x$}

\freetext(3,1){$\plus$}\freetext(6,1){$\in\lef$}

\end{graph}
\end{center}

The remaining cases follows by symmetry, and so this completes the proof.

\end{proof}

\begin{theorem}
$\hbox{Tame}_x$ is a group under the disjunctive sum.
\end{theorem}

\begin{proof}  
\noindent Case 1: Closure.  $G\plus H\in\hbox{Tame}_x$.  First assume that Left moves last on $G\plus H$, for some $G$ and $H$.  In order for Left to move last on $G\plus H$, he must either move last moving first on $G$, and moving second on $H$, or move last moving second on both $G$ and $H$.  If Left moves last moving first on $G$ and $H$, then we cannot guarantee that he will move last on $G\plus H$, since if Left makes the last move in one component, then he has to make the second move in other component.

By definition if Left moves last, moving first, on $G$, then the score will be $x$.  So if we assume that Left moved last moving first on $G$ and second on $H$, then the score from $G$ will be $x$, and 0 from $H$.  Therefore $(G\plus H)^{SL}_F=x$ if Left moves last.  By a similar argument we can show that $(G\plus H)^{SL}_F=0$ if Right moves last, and the remaining cases follow by symmetry.  Therefore the set $\tx$ is closed under the disjunctive sum.  

\noindent Case 2: Identity.  We take the identity to be all $G\in\ti$.  Since in the proof of Theorem \ref{t2}, we have already shown that if $G\in\ti$, then $G\plus X\approx X$ for all $X\in\hbox{Tame}_x$, then we have shown that there is a non-trivial identity.

\noindent Case 3: Inverse.

We will take $-G$ to be the inverse of $G$, for all $G\in\hbox{Tame}_x$.  It is enough to show that $G\plus (-G)\in\ti$ for all $G\in\hbox{Tame}_x$.  Since $\hbox{Tame}_x$ is closed, this implies that $G\plus(-G)\in\hbox{Tame}_x$.  Therefore, by definition, if the second player moves last the game is a tie.  The second player can always most last on $G\plus (-G)$ simply by playing the tweedle-dum tweedle-dee strategy.  Therefore $G\plus(-G)\in\ti$, and the theorem is proven.

\end{proof}

So what he have proven is that the games in $\hbox{Tame}_x$ do indeed form a group.  However, this is not a universe, since it is not true that if $G$ is in $\hbox{Tame}_x$, then every sub game of $G$ is also in $\hbox{Tame}_x$.  This is a significant point, but not a particularly important one for our purposes.

What we intend to show, is that if a game $G$ is in $\hbox{Tame}_x$, then we can simply use the theory of normal play combinatorial games to analyse it.  In other words, the games in $\hbox{Tame}_x$ are equivalent to normal play combinatorial games.

\begin{definition}
The function $f:\hbox{Tame}_x\rightarrow \hbox{Normal Play}$ is defined as $f(G)=\{f(G^L)|f(G^R)\}$.
\end{definition}

\begin{definition}  We define the following:
\begin{itemize}
\item{For $G,H\in\hbox{Tame}_x$, then $G\geq_x H$ if $H\plus X\in O$ then $G\plus X\in O$ for all games $X\in \hbox{Tame}_x$, where $O=L_>,R_>,L_\geq$ or $R_\geq$.}
\item{For $G,H\in\hbox{Tame}_x$, then $G\leq_x H$ if $H\plus X\in O$ then $G\plus X\in O$ for all games $X\in \hbox{Tame}_x$, where $O=L_<,R_<,L_\leq$ or $R_\leq$.}
\item{For $G,H\in\hbox{Tame}_x$ $G=_x H$, if $G\plus X\approx H\plus X$ for all $X\in\hbox{Tame}_x$.}
\end{itemize}
\end{definition}

\begin{lemma}
$G\geq_x H$ if and only if $H\leq_x G$.
\end{lemma}

\begin{proof}  The proof of this uses an identical argument as the proof of the same theorem in \cite{fs}. First let $G\geq_x H$, and let $G\plus X\in O$ for some game $X$, where $O$ is one of $L_\leq$, $R_\leq$, $L_<$ or $R_<$.  This means that $H\plus X\not\in O'$, where $O'$ is one of $L_\geq$, $R_\geq$, $L_>$ or $R_>$, since if it was this would mean that $G\plus X\in O'$, since $G\geq_x H$, therefore $H\plus X\in O$, and hence $H\leq_x G$.

A completely identical argument can be used for $H\leq_x G$, and hence $G\geq_x H$ if and only if $H\leq_x G$ and the theorem is proven.
\end{proof}

\begin{theorem}
$\hbox{Tame}_x$ is partially ordered under the disjunctive sum.
\end{theorem}

\begin{proof}  The proof of this uses an identical argument as the proof of the same theorem in \cite{fs}.  To show that we have a partially ordered set we need
3 things.

\begin{enumerate}
\item{Transitivity: If $G\geq_x H$ and $H\geq_x J$ then $G\geq_x J$.}
\item{Reflexivity: For all games $G$, $G\geq_x G$.}
\item{Antisymmetry: If $G\geq_x H$ and $H\geq_x G$ then $G=_xH$.}
\end{enumerate}

1.  Let $G\geq_xH$ and $H\geq_x J$.  $G\geq_x H$ means that if $H\plus X\in O$ this implies $G\plus X\in O$, where $O=L_\geq$, $R_\geq$, $L_>$ or $R_>$, for all games $X\in\tx$.  $H\geq_x J$, means that if $J\plus X\in O$ this implies that $H\plus X\in O$.  Since $G\geq_x H$, then this implies that $G\plus X\in O$, therefore $J\plus X\in O$ implies that $G\plus X\in O$ for all games $X$, and $G\geq_x J$.

2. Clearly $G\geq_x G$, since if $G\plus X\in O$ then $G\plus X\in O$, where  $O=L_\geq$, $R_\geq$, $L_>$ or $R_>$, for all games $X\in\tx$.

3.  First let $G\geq_x H$ and $H\geq_x G$.  $G=_xH$ means that $G\plus X\approx H\plus X$ for all $X\in\tx$.  So first let $G\plus X\in L_=$, then this implies that $H\plus X\in L_\geq$, since $H\geq_x G$.  However $H\plus X\in L_=$, since if $H\plus X\in L_>$, then this implies that $G\plus X\in L_>$, since $G\geq_x H$, therefore $G\plus X\in L_=$ if and only if $H\plus X\in L_=$.

An identical argument can be used for all remaining cases, therefore $G\plus X\approx H\plus X$ for all games $X\in\tx$, i.e. $G=_xH$.
\end{proof}

To show that we have a surjective mapping from the set $\tx$ to the set of normal play games, we need to show that the set $\tx$ can be separated into equivalence classes, in such a way that each equivalence class maps to a normal play game in canonical form.  In other words, our set $\tx$ is in fact identical to the set of normal play games.  So we begin with two definitions.

\begin{definition}
If $G\in\hbox{Tame}_x$ and $G=\{A,B,C,\dots|0|C,D,E,\dots\}$, then $A$ dominates $B$ with respect to $\hbox{Tame}_x$ if $A\geq_x B$.
\end{definition}

\begin{definition}
If $G\in\hbox{Tame}_x$ and $G=\{A,B,C,\dots|0|C,D,E,\dots\}$, then $A$ is reversible with respect to $\hbox{Tame}_x$ if $A^R\leq_x G$.
\end{definition}

\begin{theorem}\label{dom1}
Let $G, G'\in\hbox{Tame}_x$, such that $G=\{A,B,C,\dots|0|C,D,E,\dots\}$ and $G'=\{A,C,\dots|0|C,D,E,\dots\}$, then if $A\geq_x B$, $G=_x G'$.
\end{theorem}

\begin{proof}The proof of this uses an identical argument as the proof of the same theorem in \cite{fs}.   Let $G=\{A,B,C,\dots|G^S|D,E,F,\dots\}$ such that $A\geq_x B$, further let $G'=\{A,C,\dots|G^S|D,E,F,\dots\}$.  First suppose that $G\plus X\in O$, where $O=L_\geq$, $R_\geq$, $L_>$ or $R_>$, if Left moves to $B\plus X$.  This implies that $G'\plus X\in O$, since $A\geq_x B$.  Hence$G\plus X\in O$ implies that $G'\plus X\in O$, and since the Right options of $G$ and $G'$, this implies that $G'\geq_x G$.

Next suppose that $G'\plus X\in O'$ where $O'=L_\leq$, $R_\leq$, $L_<$ or $R_<$.  This implies that $G\plus X\in O'$, since the only option in $G^L$ that is not in $G'^L$ is $B$ and $B\leq_x A$, therefore $G'\leq_x G$, and $G=_x G'$.  So this means that the option $B$ may be disregarded and the proof is finished. 
\end{proof}

\begin{theorem}\label{rev1}
Let $G=\{A,B,C,\dots|G^S|D,E,F,\dots\}$, and let $A$ be reversible, with respect to $\tx$, and have
Left options of $A^R=\{W,X,Y,\dots\}$.  If
$G'=\{W,X,Y,\dots,B,C,\dots|G^S|D,E,F,\dots\}$, then $G=_xG'$.  By symmetry  if $D$ is reversible, with respect to $\tx$, and has Right options of $D^L=\{T, S, R,\dots\}$.  If $G''=\{A,B,C,\dots|G^S|T,S,R,\dots,D,E,F,\dots\}$, then $G=_xG''$.  
\end{theorem}

\begin{proof}The proof of this uses an identical argument as the proof of the same theorem in \cite{fs}.    Let $G=\{A,B,C,\dots|G^S|D,E,F,\dots\}$, where the Left options of $A^R=\{W,X,Y,\dots\}$ and let $G'=\{W,X,Y,\dots,B,C,\dots|G^S|D,E,F,\dots\}$, further let $A^R\leq G$.  If $G\plus X\in O$, where $O=L_\geq$, $R_\geq$, $L_>$ or $R_>$, when Left does not move to $A$ on $G$, then clearly $G'\plus X$ is also in $O$, since all other options for Left on $G$ are available for Left on $G'$.  

So consider the case where $G\plus X\in O$ if Left moves to $A\plus X$, then this implies that $A^R\plus X$ must also be in $O$.  This means that $G'\plus X\in O$ because $A^{RL}\subset G'^L$, and since all other options on $G'$ are the same as $G$, then $A^R\plus X\in O$ implies that $G'\plus X\in O$.  Hence if $G\plus X\in O$  then this implies that $G'\plus X\in O$, for all games $X$, i.e. $G'\geq_x G$.

Next assume that $G\plus X\in O'$, where $O'=L_\leq$, $R_\leq$, $L_<$ or $R_<$, for all games $X$.  However $A^R\leq G$, i.e. $G\plus X\in O'$ implies that $A^R\plus X\in O'$, and since $A^{RL}\subset G'^L$, and all other options on $G'$ are identical to options on $G$, this means that $G\plus X\in O'$, implies that $G'\plus X\in O'$, for all games $X$, i.e. $G'\leq_x G$.  Therefore $G=_xG'$ and the theorem is proven.
\end{proof}

\begin{theorem}
If $G,H\in\hbox{Tame}_x$ and have no dominated or reversible options with respect to $\hbox{Tame}_x$ and $G=_x H$, then $G\cong H$.
\end{theorem}

\begin{proof}  The proof of this uses an identical argument as the proof of the same theorem in \cite{fs}.  Let $G=\{A,B,C,\dots|G^S|D,E,F,\dots\}$, where the Left options of $A^R=\{W,X,Y,\dots\}$ and let $G'=\{W,X,Y,\dots,B,C,\dots|G^S|D,E,F,\dots\}$, further let $A^R\leq_x G$.  If $G\plus X\in O$, where $O=L_\geq$, $R_\geq$, $L_>$ or $R_>$, when Left does not move to $A$ on $G$, then clearly $G'\plus X$ is also in $O$, since all other options for Left on $G$ are available for Left on $G'$.  

So consider the case where $G\plus X\in O$ if Left moves to $A\plus X$, then this implies that $A^R\plus X$ must also be in $O$.  This means that $G'\plus X\in O$ because $A^{RL}\subset G'^L$, and since all other options on $G'$ are the same as $G$, then $A^R\plus X\in O$ implies that $G'\plus X\in O$.  Hence if $G\plus X\in O$  then this implies that $G'\plus X\in O$, for all games $X$, i.e. $G'\geq_x G$.

Next assume that $G\plus X\in O'$, where $O'=L_\leq$, $R_\leq$, $L_<$ or $R_<$, for all games $X$.  However $A^R\leq G$, i.e. $G\plus X\in O'$ implies that $A^R\plus X\in O'$, and since $A^{RL}\subset G'^L$, and all other options on $G'$ are identical to options on $G$, this means that $G\plus X\in O'$, implies that $G'\plus X\in O'$, for all games $X$, i.e. $G'\leq_x G$.  Therefore $G=_xG'$ and the theorem is proven.
\end{proof}

\begin{theorem}\label{map}
Let $G, H\in\hbox{Tame}_x$, if $G$ and $H$ are in canonical form with respect to $\hbox{Tame}_x$, then $G\not\cong H$ if and only if $f(G)\not\cong f(H)$, where $f(G)$ and $f(H)$ are also in canonical form.
\end{theorem}

\begin{proof}  The proof of this will be by induction.  The base case is all games of depth one or less.  So it is clear that $\{.|0|.\}\not\cong \{0\plus x|0|.\}\not\cong \{.|0|0\plus -x\}\not\cong \{0\plus x|0|0\plus -x\}$.  It should also be clear that these four games are all in canonical form, with respect to $\tx$.

For the inductive step, assume that the theorem holds for all games of depth $n$ or less, and consider two games $G$ and $H$ of depth $n+1$, such that both $G$ and $H$ are in canonical form, with respect to $\tx$, and $G\not\cong H$.

Since $G\not\cong H$, this implies that either $G^L\not\cong H^L$, or $G^R\not\cong H^R$.  Without loss of generality assume that $G^L\not\cong H^L$, since both $G^L$ and $H^L$ have depth $n$ or less, by induction this implies that $f(G^L)\not\cong f(H^L)$, which implies that $f(G)\not\cong f(H)$, and the theorem is proven.

\end{proof}

\begin{lemma}\label{sum}
For all $G,H\in\tx$, $f(G\plus H)=f(G)+f(H)$.
\end{lemma}

\begin{proof}  The proof of this will be by induction on the sum of the depths of $G$ and $H$.  The base case is trivial since $f(\{.|0.\}\plus \{.|0|.\}))=f(\{.|0|.\})+f(\{.|0|.\})=\{.|.\}$.  So for the inductive step, assume that the lemma holds for all $G$ and $H$ such that $\hbox{depth}(G)+\hbox{depth}(H)=K$, where $K$ is some positive integer.

Next consider the games $G$ and $H$ such that $\hbox{depth}(G)+\hbox{depth}(H)=K+1$. $f(G\plus H)=\{f(G^L\plus H),f(G\plus H^L)|f(G^R\plus H),f(G\plus H^R\})\}$.  However, since the sum of the depths of $G^L$ and $H$, $G$ and $H^L$, $G^R$ and $H$ and $G$ and $H^R$ is less than or equal to $K$, then by induction we have $\{f(G^L\plus H),f(G\plus H^L)|f(G^R\plus H),f(G\plus H^R\})\}=\{f(G^L)+f(H),f(G)+f(H^L)|f(G^R)+f(H),f(G)+f(H^R)\}=f(G)+f(H)$.  

So the lemma is proven.
\end{proof}

\begin{theorem}\label{t2}
For all $G,H\in\tx$, $G\geq_x H$ if and only if $f(G)\geq f(H)$.
\end{theorem}

\begin{proof}  To prove this, we first need to note one thing.  For all $G\in\tx$, $G^{SL}_F=x$ if Left moves last, and $0$ if Right moves last.  Similarly $G^{SR}_F=-x$ if Right moves last, and $0$ if Left moves last.  This should be apparent from the definition, since every move a player makes he gains precisely $x$ points, and so the final score will depend on the parity of the number of moves made, i.e. who moved last.

$G\geq_x H$ if $H\plus X\in O$ implies that $G\plus X\in O$ for all $X\in\tx$, where $O=L_\geq, R_\geq, L_>$ or $R_>$.  However, because we are in the set $\tx$, there are no games in $R_>$.  Therefore, $O=L_\geq, L_>$ or $R_=$.  $f(G)\geq f(H)$ if Left wins $f(G)+X$ whenever Left wins $f(H)+X$ for all games $X$.  Of course ``wins'' here means ``moves last''.

Our aim is to show that these two definitions are equivalent.  First let $G\geq_x H$, and consider $H\plus X\in O$.  If $O=L_>$ or $R_=$, then this implies that Left moves last on $H\plus X$, which implies that Left moves last on $G\plus X$.  By our mapping this also means that Left moves last on $f(H\plus X)$, and Left moves last on $f(G\plus X)$.  Then by Theorem \ref{map} and Lemma \ref{sum}, this is equivalent to saying that if Left wins $f(H)+X$, then Left wins $f(G)+X$ for all games $X$.

The final case to consider is if $O=L_=$.  From our definition of $\tx$, this means that Right moves last on $H\plus X$.  However, since $G\geq_x H$, then this implies that $G\plus X\in L_\geq$.  But, importantly, the converse is also true, i.e. if $G\plus X\in L_=$, then $H\plus X\in L_=$, for all $X\in\tx$.  From the previous paragraph, this is equivalent to saying that if Right wins moving first on $f(G)+X$, then Right wins moving first on $f(H)+X$ for all games $X$.

In other words, Left wins $f(G)+X$ whenever Left wins $f(H)+X$, for all games $X$, i.e. $f(G)\geq f(H)$.  Given what we know previously, if $f(G)\geq f(H)$, then this is equivalent to saying that if $H\plus X\in L_>$ or $R_=$, then $G\plus X\in L_>$ or $R_=$ for all $X\in\tx$, i.e. $G\geq_x H$.

Therefore $G\geq_x H$, if and only if $f(G)\geq f(H)$, and the theorem is proven.
\end{proof}

So what this sequence of theorems and proofs have shown is that the set $\tx$ is completely identical to the set of normal play games.  That is, we can use scoring play theory on normal play games, and normal play theory on games in $\tx$, and we will get completely identical results.

\subsection{Surreal Numbers}

Given that we now have a surjective mapping from a subset of scoring games to the set of normal play games, we can now re-build the surreal numbers using the definitions from scoring play combinatorial game theory.

\begin{definition}
A game $Y=\{Y^L\plus x|0|Y^R\plus -x\}\in\hbox{Tame}_x$ is a number if for all $y^L\in Y^L$ and $y^R\in Y^R$, $y^L\not\geq_x y^R$.
\end{definition}

Using this definition it is possible to redefine the surreal numbers using the definitions from scoring play game theorem.  For instance the game born on day 0 is simply $\{.|0|.\}=0$, then we have $\{0\plus x|0|.\}=1, \{.|0|0\plus -x\}=-1$ born on day 1.

We can continue in this manner, and we simply get the set of surreal numbers again.  The purpose of this is not to propose an improvement over the Conway number system.  Indeed it is clear that the original Conway numbers are much more convenient to use on a notational level.  We simply intended to show that one can use scoring play theory to derive the surreal numbers.

\subsection{Mapping to Mis\`ere Play}

The mapping to mis\`ere play is almost identical to the mapping to normal play.  The only difference is the set that we use to do it.

\begin{definition}\label{tame}
$\hbox{Tame}_{-x}=\{G|G=\{G^L\plus -x|0|G^R\plus x\}, G^L,G^R\in\hbox{Tame}_{-x}\}$
\end{definition}

The important difference to note here is that $G^{SL}_F=0$ if Right moves last, and $-x$ if Left moves last.  Likewise, $G^{SR}_F=0$ if Left moves last, and $x$ if Right moves last.  This means that both players are trying \emph{not} to move last.  We will now show that this set is identical to the set of mis\`ere games.

\begin{theorem}
For all $G\in\tm$ if $G\not\cong 0$, then $G\neq_{-x}0$.
\end{theorem}

\begin{proof}  The proof of this is almost identical to the proof found to a similar theorem in \cite{gamo}.  Consider the game $G$.  This game is an $\lef$ position, next consider a game $H$, and without loss of generality assume that $H^L\not\emptyset$.  Right moving first on $G\plus H$, will move to $G^R\plus H$.  Left then must move to $G^R\plus H^L$, Right can then guarantee to at least tie moving second by moving to $G^RR\plus H^L$.  Therefore $G\plus H$ is not an $\lef$ position, i.e. $G\plus H\not\approx G$, and the theorem is proven.

\begin{center}
\begin{graph}(2,4)

\roundnode{1}(0,0)\roundnode{2}(0.5,.5)\roundnode{3}(1.5,1.5)
\roundnode{4}(2,2)\roundnode{5}(1.5,2.5)\roundnode{6}(1,3)

\edge{1}{2}\edge{2}{3}[\graphlinedash{1 4}]\edge{3}{4}
\edge{4}{5}\edge{5}{6}

\freetext(0,1.5){$G=$}

\nodetext{1}{$(2-n)x$}\nodetext{2}{$(3-n)x$}\nodetext{3}{$x$}
\nodetext{4}{$2x$}\nodetext{5}{$x$}\nodetext{6}{0}

\end{graph}
\end{center}

\end{proof}

\begin{definition}
The function $g:\hbox{Tame}_{-x}\rightarrow \hbox{Mis\`ere Play}$ is defined as $g(G)=\{g(G^L)|g(G^R)\}$.
\end{definition}

\begin{theorem}
For any outcome classes $\mathcal{X},\mathcal{Y}$ and $\mathcal{Z}$, there is a $G\in\mathcal{X}$, $H\in\mathcal{Y}$ such that $G\plus H\in\mathcal{Z}$.
\end{theorem}

\begin{proof}  The proof of this uses exactly the same examples that can be found the appendix of \cite{gamo}, with the vertices of the game trees labelled appropriately.  So we will not reproduce it here.
\end{proof}

The entire same procedure, and theorems, from Section \ref{normal} can be used to derive canonical forms.  So again, we will not reproduce it here.  Instead we will simply state the following theorem.

\begin{theorem}
Let $G, H\in\hbox{Tame}_{-x}$, if $G$ and $H$ are in canonical form with respect to $\hbox{Tame}_{-x}$, then $G\not\cong H$ if and only if $f(G)\not\cong f(H)$, where $f(G)$ and $f(H)$ are also in canonical form.
\end{theorem}

\begin{proof}  The proof of this is almost identical to the proof of Theorem \ref{map}.
\end{proof}

\begin{theorem}
For all $G,H\in\tm$, $G\geq_{-x}H$ if and only if $g(G)\geq g(H)$.
\end{theorem}

\begin{proof}  The proof of this is identical to Theorem \ref{t2}.
\end{proof}

Again, what we have shown is that we can use mis\`ere play theory to analyse games in $\tm$, but more importantly, scoring play theory to analyse mis\`ere games.

\section{Conclusion}

With this paper, what we have shown is that all combinatorial games, normal play, mis\`ere play and scoring play can be analysed using scoring play combinatorial game theory.  This means that scoring play combinatorial game theory is in fact, a complete theory for combinatorial games.  We hope that this inspires further research into the area.

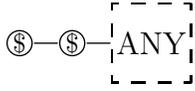
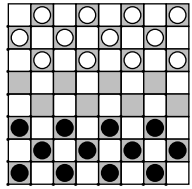
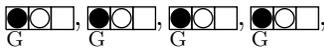
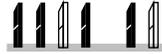
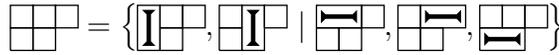
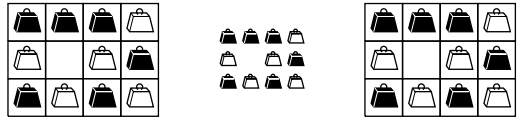
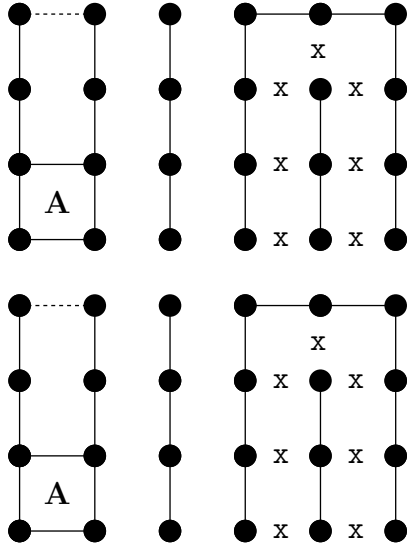
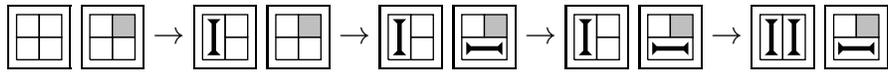
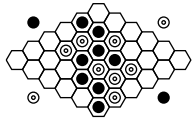
\end{document}